\documentclass[a4paper,11pt]{article}
\usepackage[french]{babel}
\usepackage[latin1]{inputenc}

\title{\textbf{Sur la rigidité de polyèdres hyperboliques en dimension $3$ : cas de volume fini, cas hyperidéal, cas fuchsien. }}
\author{Mathias Rousset}
\date{avril-juin $2002$}
\begin{document}
\newtheorem{The}{Théorème}
\newtheorem{Pro}{Proposition}
\newtheorem{Lem}{Lemme}
\newtheorem{Cor}{Corollaire}
\newtheorem{Def}{Définition}
\newtheorem{Prodef}{Proposition-définition}

\maketitle
\paragraph{Abstract.} A hyperbolic semi-ideal
polyedron is a polyedron whose vertices lie inside the hyperbolic
space $\mathbf{H}^{3}$ or at infinity.\\
A hyperideal polyedron is, in the projective model, the
intersection of $\mathbf{H}^{3}$ with a projective polyhedron
whose vertices all lie outside of $\mathbf{H}^{3}$, and whose
edges all meet $\mathbf{H}^{3}$.\\
We classify semi-ideal polyhedra in terms of their dual metric,
using the results of Rivin and Hodgson in \cite{comp} et
\cite{idea}. This result is used to obtain the classification of
hyperideal polyhedra in terms of their combinatorial type and
their dihedral angles. These two results are generalized to the
case of fuchsian polyhedra.
\paragraph{Résumé.}Un polyèdre hyperbolique semi-idéal est un
polyèdre dont les sommets sont dans l'espace hyperbolique
$\mathbf{H}^{3}$ ou à l'infini.\\
Un polyèdre hyperbolique hyperidéal est, dans le modèle projectif,
l'intersection de $\mathbf{H}^{3}$ avec un polyèdre projectif dont
les sommets sont tous en dehors de $\mathbf{H}^3$, et dont toutes
les arêtes rencontrent $\mathbf{H}^{3}$.\\
Nous classifions les polyèdres semi-idéaux en fonction de leur
métrique duale, d'après les résultats de Rivin et Hodgson dans
\cite{comp} et \cite{idea}. Nous utilisons ce résultat pour
retrouver la classification des polyèdres hyperidéaux en terme de
leurs angles dièdres. Nous généralisons ces résultats au cas des
polyèdres fuchsiens.
\newpage
\tableofcontents
\newpage
\section{Introduction}
L'objet principal d'étude de cet article est la rigidité des
polyèdres hyperboliques. Cauchy a montré la rigidité globale des
polyèdres euclidiens, les faces étant fixées à isométrie prés.
D'une manière générale, un problème de rigidité de polyèdres peut
se voir de la manière suivante:
\begin{itemize}
\item On considère un certain espace de polyèdres convexes, en général on
fixe le nombre de plans formant le polyèdre, ou bien (et c'est
très différent) la combinatoire du polyèdre (ie la décomposition
en arêtes, sommets et faces).
\item on considère une application de ces polyèdres vers certaines
de leurs caractéristiques, typiquement la longueur des arêtes, la
valeur des angles dièdres ou la métrique induite sur le polyèdre,
etc...
\item on se demande si l'application est localement injective (rigidité locale), ou globalement
injective (rigidité globale).
\end{itemize}
Dans le cas de Cauchy, on fixe la combinatoire du polyèdre, et on
regarde la longueur des arêtes et la valeur des angles des faces.

On peut être plus ambitieux et essayer de caractériser un ensemble
de polyèdres, c'est à dire le classifier par certaines de ses
caractéristiques. On peut par exemple essayer de mettre en
bijection les polyèdres de combinatoire fixé $\Gamma$ avec
certaines valeurs de leurs angles dièdres. Ou bien encore, on peut
essayer de mettre en bijection les polyèdres avec $n$ sommets avec
un certain ensemble de métriques sur la sphère avec $n$
singularités représentant leur géométrie intrinsèque.

\subsection{Résultats principaux}
Commençons par rappeler les résultats importants sur la rigidité
des polyèdres. Le premier est du à Cauchy et s'énonce ainsi :
\begin{The}[Cauchy]
Soit $\mathcal{P}_{euc,\Gamma}$ l'ensemble des polyèdres
euclidiens de combinatoire $\Gamma$. L'application qui à un
polyèdre de $\mathcal{P}_{euc,\Gamma}$ associe la longueur de ses
arêtes et la valeur des angles de ses faces est globalement
injective.
\end{The}
Andreev a montré la rigidité des polyèdres hyperboliques de volume
fini dans \cite{andreev1} pour le cas compact, puis dans
\cite{andreev2} dans le cas de volume fini :
\begin{The}[Andreev]
Soit $\mathcal{P}_{and,\Gamma}$ l'ensemble des polyèdres
hyperboliques de volume fini de combinatoire $\Gamma$, et ayant
des angles dièdres inférieurs à $\pi / 2$. L'application qui à un
polyèdre de $\mathcal{P}_{and,\Gamma}$ associe la valeur de ses
angles dièdres est globalement injective.
\end{The}
Andreev donne un ensemble de conditions caractérisant les valeurs
d'angles dièdres atteints.

On a ensuite une caractérisation très générale des polyèdres
hyperboliques de volume fini dont \emph{le nombre de faces est
fixé mais pas la combinatoire}; caractérisation en fonction de la
métrique duale (la métrique du "polyèdre dual" plongée dans
l'espace de Sitter ou de manière équivalente, la troisième forme
fondamentale du polyèdre). On notera $\mathcal{P}_{n}$ l'ensemble
des polyèdres hyperboliques compacts ayant $n$ faces, et
$\mathcal{P}^{si}_{n}$ l'ensemble des polyèdres hyperboliques de
volume fini (on dira aussi semi-idéal)ayant $n$ faces. On définit
maintenant les espaces métriques qui vont caractériser ces
polyèdres :
\begin{Def} On notera $\mathcal{S}_{n}$ l'espace des espaces
métriques $h$, tel que h soit homéomorphe à la sphère
$\mathbf{S}^{2}$, de courbure constante $1$ partout à l'exception
d'un nombre fini de points où il présente une singularité conique.
Les singularités sont numérotées, et les espaces métriques sont
définis aux isométries préservant les singularités prés.
\end{Def}
\begin{Def} On notera $\mathcal{M}_{n}$ le sous-ensemble de
$\mathcal{S}_{n}$ des métriques telles que :
\begin{itemize}
\item Toutes les singularités ont une courbure strictement
négative.
\item Toute géodésique fermée a une longueur strictement supérieure à $2\pi$.
\end{itemize}
On notera $\mathcal{M}^{si}_{n}$ le sous-ensemble de
$\mathcal{S}_{n}$ des métriques telles que:
\begin{itemize}
\item Toutes les singularités ont une courbure strictement
négative.
\item toute géodésique fermée a une longueur
supérieure ou égale à $2\pi$.
\item le nombre de géodésiques fermées de longueur $2\pi$ est fini
et chacune d'entre elles sépare l'espace en deux composantes
connexes dont l'une au moins est isométrique à un hémisphère.
\end{itemize}
\end{Def}
A tout polyèdre compact $P$ de $\mathcal{P}_{n}$ on associe son
dual $P^{*}$ muni de sa métrique induite (cf. section 3). En tant
qu'objet intrinsèque, c'est un élément de $\mathcal{S}_{n}$. On
note $\Phi$ l'application ainsi définie. On est en mesure de citer
des résultats de Rivin et Hodgson \cite{comp} et Rivin
\cite{idea}, respectivement :
\begin{The}[Rivin et Hodgson]\label{Rivi}
L'application $\Phi$ ainsi définie est un homéomorphisme de
$\mathcal{P}_{n}$ sur $\mathcal{M}_{n}$.\newline Autrement dit:
Les polyèdres compacts sont en bijection avec les métriques duales
admissibles.(Existence et unicité).
\end{The}

\begin{The}[Rivin]\label{Rivi}
Le prolongement par continuité de $\Phi$ à
$\tilde{\Phi}:\mathcal{P}^{si}_{n}\rightarrow \mathcal{S}_{n}$ a
exactement pour image $\mathcal{M}^{si}_{n} $.\newline Autrement
dit: les polyèdres semi-idéaux ont exactement pour métriques
duales les métriques admissibles de volume fini.(Existence).
\end{The}

Antérieurement, Alexandrov \cite{alex} a donné une caractérisation
similaire pour les polyèdres compacts ayant n sommet, simplement
en fonction de leur métrique induite; la condition de borne
inférieure pour les longueurs des géodésiques fermées n'apparaît
alors pas.

Enfin, citons le résultat de Bonahon et Bao \label{hype}, qui
permet de caractériser les polyèdres hyperboliques hyperidéaux
(avec leurs sommets à l'infini et "au-delà" de l'infini, cf.
section 4) de combinatoire donnée $\Gamma$ (dont l'ensemble est
noté $\mathcal{P}_{\Gamma}$) par la donnée de leurs angles dièdres
:
\begin{The} [Bonahon et Bao] \label{bonha}
On a existence et unicité des polyèdres hyperidéaux (on autorise
des sommets idéaux) de combinatoire donnée $\Gamma$ dont
les angles dièdres extérieurs $\theta_{e_{i}}$ ($\pi $ moins l'angle dièdre) vérifient les deux conditions suivantes :\\
(\emph{C1})$\sum_{i=1}^{n}\theta_{e_{i}}\geq 2\pi$ pour toute
courbe fermée $\gamma$ plongée dans le graphe dual $\Gamma^{*}$ et
passant par les arêtes $e_{1},...,e_{n}$; l'égalité ayant lieu
très exactement quand $\gamma$ est le bord d'une face du graphe
dual $\Gamma^{*}$ correspondant à un sommet idéal du polyèdre initial.\\
(\emph{C2}) $\sum_{i=1}^{n}\theta_{e_{i}} > \pi$ pour tout chemin
$\gamma$ plongé dans $\Gamma^{*}$ et joignant deux sommets de
$\Gamma^{*}$ ayant une même face adjacente, mais tel que $\gamma$
ne soit pas tout entier contenu dans le bord de cette face.\\
\end{The}

Cet article contient les deux résultats principaux suivants :
\begin{itemize}
\item La caractérisation des polyèdres convexes hyperboliques à $n$ faces semi-idéaux
(avec des sommets à l'infini) en fonction de leur métrique duale.
Cette caractérisation comporte comme cas particulier celle des
polyèdres compacts et des polyèdres idéaux.
\item Une nouvelle démonstration de la caractérisation des polyèdres convexes hyperboliques de combinatoire $\Gamma$ hyperidéaux (les sommets
sont à l'infini ou au-delà) par leurs angles dièdres.
\end{itemize}
L'idée est de déduire le premier résultat des idées de Rivin et
Hodgson \cite{comp} et \cite{idea}(section 3), l'apport nouveau du
article est l'obtention de la rigidité dans le cas général
semi-idéal (section 3), en utilisant la transformation de
Pogorelov (section 2).\\
En "tronquant" les sommets hyperidéaux, on en déduit le second
énoncé (section 5). Cette partie constitue une preuve indirecte
simple, du résultat démontré en détail dans la prépublication de
Bonahon et Bao \cite{hype}.\\
On montrera que ces caractérisations se généralisent très
exactement au cas des polyèdres hyperboliques fuchsiens de genre
$g \geq 2$ (cf. section 5). En particulier, on obtient ainsi le
résultat nouveau suivant de caractérisation de polyèdres hyper
idéaux fuchsiens :
\begin{The} \label{fuchshyp}
On a existence et unicité des polyèdres hyperidéaux fuchsiens de
genre $g \geq 2$ (on autorise des sommets idéaux) de combinatoire
donnée $\Gamma$ dont
les angles dièdres extérieurs $\theta_{e_{i}}$ ($\pi $ moins l'angle dièdre) vérifient les deux conditions suivantes :\\
(\emph{C1})$\sum_{i=1}^{n}\theta_{e_{i}}\geq 2\pi$ pour toute
courbe fermée \emph{contractile}$\gamma$ plongée dans le graphe
dual $\Gamma^{*}$ et passant par les arêtes $e_{1},...,e_{n}$;
l'égalité ayant lieu très exactement quand $\gamma$ est le bord
d'une face du graphe
dual $\Gamma^{*}$ correspondant à un sommet idéal du polyèdre initial.\\
(\emph{C2}) $\sum_{i=1}^{n}\theta_{e_{i}} > \pi$ pour tout chemin
$\gamma$ plongé dans $\Gamma^{*}$ et joignant deux sommets A et B
de $\Gamma^{*}$ ayant une même face F adjacente, mais tel que
$\gamma$ ne soit pas tout entier contenu dans le bord de cette
face. On demande que $\gamma$
soit homotope au chemin bordant F et joignant A à B.\\
\end{The}
L'analogue des résultats de Rivin et Hodgson \cite{comp} et
\cite{idea} dans le cas fuchsien étant faits par Schlenker dans
\cite{fuchs}.
\section{La transformation de Pogorelov}
Cette section est consacrée à l'étude de la transformation de
Pogorelov. Cette transformation est caractérisée par la propriété
suivante : On note $\mathbf {R}^{n+1}_{p}$ l'espace de Minkowski
de dimension $n+1$ et de signature $(n+1-p,p)$, et $\Omega_{0}$
l'intersection d'un demi-espace de vecteur normal de type
strictement espace ou strictement temps avec une des deux sphères
unités ($p \neq 0$) .\newline La transformation de Pogorelov prend
deux objets de $\Omega_{0}$ ayant la même géomètrie intrinsèque
induite et les plongent dans $\mathbf{R}^{n}_{p}$ ou
$\mathbf{R}^{n}_{p-1}$, les images ayant alors la même géométrie
intrinsèque induite. De plus, les images sont identifiables à
isométrie prés si et seulement si les objets initiaux le sont.
Elle va nous permettre de ramener des questions de rigidité de
polyèdres dans l'espace de Sitter, au même problème dans l'espace
euclidien.
\subsection{Espaces à courbure constante}\label{cons}
On se propose dans cette sous-section de rappeler succinctement
quelques propriétés des espaces homogènes riemanniens et
lorentziens à courbure constante. L'espace ambiant est l'espace de
Minkowsky $\mathbf{R}^{n+1}_{p}$ muni de sa métrique canonique
:$$\parallel x
\parallel^{2}=-x_{1}^{2}-...-x_{p}^{2}+x_{p+1}^{2}+...+x_{n+1}^{2}$$
On note $\mathbf{S}^{n}_{p}=\{x\in \mathbf{R}^{n+1}_{p},\parallel
x
\parallel = 1\} $ la sphère unité dite "à courbure positive" et
$\mathbf{H}^{n}_{p-1}=\{x\in \mathbf{R}^{n+1}_{p},\parallel x
\parallel = -1\} $ la sphère unité dite "à courbure négative", toutes deux munies de leur métrique induite.\newline On
énonce sans démonstration :
\begin{The}\label{unit} La sphère $\mathbf{S}^{n}_{p}$ est une variété
pseudo-riemannienne complète de dimension $n$ et de signature
$(n-p,p)$. Elle est homogène et à pour courbure constante $1$. Son
groupe d'isométrie est $O(n+1-p,p)$, le groupe des isométries
linéaires de $\mathbf{R}^{n+1}_{p}$ et son groupe d'isotropie est
$O(n+1-p,p-1)$.\newline La sphère $\mathbf{H}^{n}_{p-1}$ est une
variété pseudo-riemannienne complète de dimension $n$ et de
signature (n+1-p,p-1). Elle est homogène et à pour courbure
constante $-1$. Son groupe d'isométrie est $O(n+1-p,p)$ et son
groupe d'isotropie est $O(n-p,p)$
\end{The}

\subsection{La transformation projective} On se donne un demi-espace de
$\mathbf{R}^{n+1}_{p}$ de normale $x_{0}$ avec $\parallel x_{0}
\parallel^{2}= \epsilon = \pm 1$. On note $\Omega_{0}$ l'ensemble des
points de norme carrée $\mu = \pm 1$ de ce demi-espace. Donc
$\Omega_{0}=\{x \in \mathbf{S}^{n}_{p},sgn(<x,x_{0}>)=\epsilon \}$
si $\mu$ est positif, $\Omega_{0}=\{x \in
\mathbf{H}^{n}_{p-1},sgn(<x,x_{0}>)=\epsilon\}$ dans le cas
contraire.

On considère l'hyperplan affine passant par le point $x_{0}$ et
orthogonal au vecteur $x_{0}$. C'est un modèle affine pour
l'espace projectif $\mathbf{RP}^{n}$ des droites vectorielles de
$\mathbf{R}^{n+1}_{p}$. Chaque droite de cet ensemble intersecte
$\Omega_{0}$ en au plus un point. Cette construction fournit
l'application projective $\phi$ de $\Omega_{0}$ vers
$\mathbf{R}^{n}\subset \mathbf{RP}^{n}$.

On énonce alors sans démonstration :
\begin{Pro}[transformation projective]\label{proj}Avec les notations précédentes, la
transformation projective $\phi:\Omega_{0}\rightarrow
\mathbf{R}^{n}$ est un difféomorphisme de $\Omega_{0}$ sur son
image. On appelle $x_{0}$ le "centre" de la projection. Si on
considère $\phi(x)$ comme vecteur de $ \mathbf{R}^{n}_{p}$ il a
pour expression explicite :
$$\phi(x)=\frac{\epsilon x-<x,x_{0}>x_{0}}{<x,x_{0}>}$$
Enfin, si $\rho$ est une isomètrie de $\Omega_{0}$ laissant le
centre invariant, $\phi \circ \rho \circ \phi^{-1}$ est une
isométrie de $\mathbf{R}^{n}_{p}$ laissant le centre invariant.
\end{Pro}
\paragraph{Remarque.}On a la formule suivante pour $\phi^{-1}$ :
$$\phi^{-1}(y)=(y+x_{0})(\mu \parallel y+x_{0} \parallel
^{2})^{\frac{1}{2}}$$.

\subsection{La transformation de Pogorelov globale}
On conserve les notations précédentes.
\begin{The}[Application de Pogorelov globale]Il existe une
application $\Phi$ : $\Omega_{0} \times \Omega_{0} \rightarrow
\mathbf{R}^{n}_{p} \times \mathbf{R}^{n}_{p}$, dont la formule
explicite est :
$$\Phi(x,x')=2(\frac{\epsilon
x-<x,x_{0}>x_{0}}{<x+x',x_{0}>},\frac{\epsilon
x'-<x',x_{0}>x_{0}}{<x+x',x_{0}>}$$ et qui vérifie les deux
propriétés suivantes :
\begin{enumerate}
\item Soit $V \in T_{(x,x')}(\Omega_{0} \times \Omega_{0})$. On a
que $\parallel \pi_{1}V \parallel^{2}=\parallel \pi_{2}V
\parallel^{2}$, si et seulement si $\parallel
\pi_{1}(T_{(x,x')}\Phi(V))
\parallel^{2}=\parallel \pi_{2}(T_{(x,x')}\Phi(V)\parallel^{2}$; $\pi_{1}$
et $\pi_{2}$ étant respectivement la première et seconde
projection d'un produit cartésien.
\item $\Phi$ commute aux isométries laissant le centre $x_{0}$
invariant, au sens où si $\rho$ appartient au groupe d'isotropie
de $x_{0}$, on a:
$$\Phi \circ ( \rho \times \rho )=(\phi \circ \rho \circ
\phi^{-1}) \times (\phi \circ \rho \circ \phi^{-1})) \circ \Phi$$
$\phi \circ \rho  \circ \phi^{-1}$ étant une isométrie.
\end{enumerate}

\end{The}
\paragraph{Remarque.}$\Phi_{\mid \Delta}=(\phi \times \phi)_{\mid
\Delta}$ , où $\Delta$ est la diagonale de $\Omega_{0} \times
\Omega_{0}$.

\paragraph{\emph{Démonstration.}}
\begin{enumerate}
\item On se donne $(X,X') \in T_{(x,x')}(\Omega_{0} \times
\Omega_{0})$ et on calcul la différentielle de $\Phi$
directement:\begin{eqnarray*}&&\frac{1}{2}<x+x',x_{0}>^{2} \pi_{1}
(T_{(x,x')}
\Phi(X,X'))=\\
 &&\qquad\epsilon<x+x',x_{0}>X-\epsilon<X+X',x_{0}>x\\
 &&\qquad+(<x,x_{0}><X',x_{0}> -<x',x_{0}><X,x_{0}>)x_{0}
\end{eqnarray*}
Puis sa norme au carré, en utilisant $<X,x>=0$.
\begin{eqnarray*}&&\frac{1}{4}<x+x',x_{0}>^{4}\parallel \pi_{1}
(T_{(x,x')} \Phi(X,X'))\parallel^{2}=\\
 &&\qquad<x+x',x_{0}>^{2}\parallel X \parallel^{2} + \mu<X+X',x_{0}>^{2} \\
 &&\qquad-\epsilon(<x,x_{0}><X',x_{0}> -<x',x_{0}><X,x_{0}>)^{2}
\end{eqnarray*}
Enfin on obtient la différence des normes des deux projections :
\begin{displaymath}
\frac{1}{4}<x+x',x_{0}>^{2}(\parallel \pi_{1} (T
\Phi(X,X'))\parallel^{2}-\parallel \pi_{2}
(T\Phi(X,X'))\parallel^{2})=\parallel X \parallel^{2}-\parallel X'
\parallel^{2}
\end{displaymath}
Ce qui suffit à compléter la preuve car $<x+x',x_{0}>$ n'est
jamais nul.
\item L'égalité s'obtient en effectuant un calcul explicite et en
exploitant le linéarité de $\rho$, le fait qu'il laisse $x_{0}$
invariant, et le fait qu'il conserve le produit scalaire.
\end{enumerate}
La démonstration est ainsi achevée.$\diamond$ \\
Voici un corollaire qui apporte l'interprétation géométrique de la
transformation :
\begin{Cor}
Soit $F$ et $F'$ deux sous-variétés difféomorphes de $\Omega_{0}$,
plongements d'une variété $S$ par $\phi$ et $\phi '$
respectivement, munis de la métrique induite, et $\tilde{F}$ et
$\tilde{F'}$ le couple de sous-variétés image du couple $(F,F')$
par la transformation de Pogorelov, munis de la métrique induite
par $\mathbf{R}_{p}^{n}$. On a les deux propriétés suivantes :
\begin{enumerate}
\item $F$ et $F'$ ont la même géométrie intrinsèque si et
seulement si $\tilde{F}$ et $\tilde{F'}$ aussi.
\item $F$ et $F'$ sont identiques à isométrie globale préservant le centre de $\Omega_{0}$ prés, si et
seulement si $\tilde{F}$ et $\tilde{F'}$ aussi.
\end{enumerate}
\end{Cor}
\paragraph{Démonstration.}
\begin{enumerate}
\item Soit $S$ tel que $F$ et $F'$ soient deux plongements de $S$ dans
$\Omega_{0}$ par $\phi$ et $\phi'$ respectivement. $\tilde{F}$ et
$\tilde{F'}$ sont deux plongements de S dans $\mathbf{R}_{p}^{n}$
par $\pi_{1}\Phi \circ(\phi, \phi')$ et $\pi_{2}\Phi \circ(\phi,
\phi')$ respectivement. La première propiété de la transformation
dit exactement que les premiers plongements sont isométriques si
et seulement si les seconds le sont.
\item D'après la seconde propriété, si $F'$ est l'image de $F$ par l'isométrie $\rho$,
$\tilde{F'}$ est l'image de $\tilde{F}$ par $\phi \circ \rho \circ
\phi^{-1}$. On obtient la réciproque en conjuguant par $\phi^{-1}$
l'isométrie de $\mathbf{R}_{p}^{n}$  préservant le centre.
\end{enumerate}
Ceci achève la démonstration du corollaire.$\diamond$
\subsection{La transformation de Pogorelov infinitésimale}
\begin{The}[Application de Pogorelov infinitésimale.]\label{Pogi}Il existe un
morphisme de fibrés : $\tilde{\phi}:T\Omega_{0}\rightarrow
T\mathbf{R}_{p}^{n}$ de base $\phi$, s'écrivant explicitement avec
les notations précédentes :
\begin{eqnarray*}
\tilde{\phi}:  T \Omega_{0} & \rightarrow & T \mathbf{R}_{p}^{n} \\
(x,v) & \mapsto & (\phi(x)=\frac{\epsilon
x-<x,x_{0}>x_{0}}{<x,x_{0}>},\frac{\epsilon
v-<v,x_{0}>x_{0}}{<x,x_{0}>})
\end{eqnarray*}
et vérifiant les propriétés suivantes :
\begin{enumerate}
\item Si $F$ est une sous-variété de $\Omega_{0}$ muni d'un champ
de vecteur $V$, $V$ est une  déformation isométrique (ie de
dérivée de Lie nulle sur $F$),si et seulement si $\tilde{\phi}(V)$
est une déformation isométrique de $\phi(F)$.
\item $V$ est la restriction sur $F$ d'une isométrie
infinitésimale (ie un champ de Killing) $\xi$, ssi $
\tilde{\phi}(V)$ est la restriction à $\phi(F)$ de l'isométrie
infinitésimale $\tilde{\phi}(\xi)$.
\end{enumerate}
\end{The}
\paragraph{Démonstration.}On peut donner deux démonstrations du
théorème. La première en dérivant la formule de Pogorelov globale,
la seconde par le calcul direct. On rappelle que la dérivée de Lie
est la deux forme symétrique
$\mathcal{L}(Y,Z)=<D_{Y}V,Z>+<D_{Z}V,Y>$
\\\emph{Première démonstration.}On se donne $f_{t}:
S\hookrightarrow \Omega_{0}$ une famille différentiable de
plongements pour $t\in ]-\epsilon,\epsilon[$, isométriques entre
eux ( la deux-forme $< Tf_{t},Tf_{t}>$ ne dépend pas de t), et
telle que $F\equiv f_{0}$ et
$V(f_{0}(p))=\frac{d}{dt}_{|t=0}(f_{t}(p))$, p étant un point
courant de $S$. La composition $\Phi(f_{0},f_{t})$, donne d' après
le théorème de Pogorelov, une famille de couples isométriques de
plongements de $S$ dans $\mathbf{R}_{p}^{n}$. Soient $Y$ et $Z$
deux vecteurs de l'espace tangent en $p$ à $S$, et soient
($Y_{t}^{1}$, $Y_{t}^{2})$) et ($Z_{t}^{1}$, $Z_{t}^{2})$) leurs
images par $\Phi(f_{0},f_{t})$($Y_{0} \equiv Y_{0}^{i})$.
Progorelov nous dit très exactement que :$<Y_{t}^{1},
Z_{t}^{2}>$=$<Y_{t}^{2}, Z_{t}^{2}>$ pour tout t. Dérivons cette
égalité, il vient
:$$<T_{p}V_{1}(Y_{0}),Z_{0}>+<T_{p}V_{1}(Z_{0}),Y_{0}>=<T_{p}V_{2}(Y_{0}),Z_{0}>+<T_{p}V_{2}(Z_{0}),Y_{0}>$$
où on a posé :$V_{i}=\pi_{i}T_{p}\Phi(O,V)$. Cette formule dit
très exactement que $V_{1}-V_{2}$ à une dérivée de Lie nulle sur
$\phi(F)$. On pose donc $\tilde{\phi}(V)=V_{1}-V_{2}$ et un calcul
explicite donne la formule.$\diamond$\\
\emph{Deuxième démonstration.} On se donne un point $p$ de F et
$y_{t}$ un chemin de $F$ tel que $y_{0}=p$. On note
$Y=\frac{dy_{t}}{dt}_{|t=0}$. On cherche à calculer la dérivée de
Lie (une deux forme symétrique) de $\tilde{\phi}(V)$ sur deux fois
$T\phi(Y)=\frac{d\phi(y_{t})}{dt}=$ en fonction de la dérivée de
Lie de $V$ sur deux fois $Y$. Un calcul explicite donne en $t=0$
la
formule:$$<\frac{d}{dt}\tilde{\phi}(V(y_{t})),\frac{d\phi(y_{t})}{dt}>=<\tilde{\phi}(\frac{dV(y_{t})}{dt})
-\tilde{\phi}(V)\frac{<Y,x_{0}>}{<p,x_{0}>},\tilde{\phi}(Y)-\phi(p)\frac{<Y,x_{0}>}{<p,x_{0}>}>$$
Puis en développant et en utilisant les identités
\mbox{$<V,p>=0$,$<Y,p>=0$},\\\mbox{$<\frac{dV(y_{t})}{dt},p>+<V,Y>=0$},
il vient :
$$<\frac{d}{dt}\tilde{\phi}(V(y_{t})),\frac{d\phi(y_{t})}{dt}>=<p,x_{0}>^{-2}<\frac{dV(y_{t})}{dt},Y>$$
Ce qui constitue la relation souhaitée et prouve que
$\tilde{\phi}(V)$ est isométrique si et seulement si $V$ l'est.\\
La deuxième propriété découle immédiatement de la seconde en
faisant $F=\Omega_{0}$. $\diamond$
\section{Une caractérisation des polyèdres hyperboliques convexes de volume
fini} Le but de cette section est de résumer les résultats de Igor
Rivin et Craig Hodgson contenus dans \cite{comp} et \cite{idea},
puis d'ajouter un argument d'unicité grâce à la transformation de
Pogorelov.
\subsection{Polyèdres}
On se place dans l'espace hyperbolique $\mathbf{H}^{3}$ de
dimension 3.
\begin{Def}[Polyèdres convexes]
On appelle \emph{polyèdre hyperbolique convexe compact }
l'intersection compacte d'une collection de n demi-espaces de
$\mathbf{H}^{3}$. Tous ses sommets sont situés dans
$\mathbf{H}^{3}$ \\ On appelle \emph{polyèdre hyperbolique convexe
de volume fini ou semi-idéal } l'intersection de volume fini d'une
collection de n demi-espaces de $\mathbf{H}^{3}$. Ses sommets sont
situés dans $\mathbf{H}^{3}$ où à l'infini. On appelle
\emph{polyèdre hyperbolique convexe idéal}, un polyèdre convexe de
volume fini dont tous les sommets sont situés à l'infini.
\end{Def}
\begin{Def}
On appelle \emph{combinatoire} du polyèdre $P$ la décomposition
cellulaire du polyèdre selon ses arêtes, ses sommets et ses faces.
On appelle \emph{combinatoire duale} de $\Gamma$, la combinatoire
dont les faces correspondent aux sommets de $\Gamma$, les arêtes
aux arêtes, et les sommets aux faces. (Ex: la combinatoire duale
du cube est celle de l'octaèdre).
\end{Def}

\paragraph{Remarque.}On définit de même les polyèdres convexes
sphériques. On a aussi une notion de polyèdres convexe dans
l'espace de se Sitter $\mathbf{S}^{3}_{1}$ :
\begin{Def}On appelle \emph{polyèdre convexe compact } dans
l'espace de Sitter l'intersection compacte d'une collection de n
demi-espaces définis par des plans orientés de type-espace. La
métrique induite sur ces polyèdres est donc riemannienne.
\end{Def}
Une manière simple de visualiser topologiquement ces polyèdres est
la suivante, en notant $can_{\mathbf{S}^{2}}$ la métrique
naturelle sur $\mathbf{S}^{2}$ :

\begin{Pro}On considère l'espace de Sitter comme la variété
différentielle $\mathbf{S}^{2}\times \mathbf{R}$ muni de
$(-dt^{2}+ch^{2}t\,can_{\mathbf{S}^{2}})$, et on note $p$ la
projection sur la sphère. Soit $P$ un polyèdre compact, alors $p$
constitue un homéomorphisme (aux singularités prés) de P sur
$S^{2}$.
\end{Pro}
\subsection{Dualité}
On va définir une application notée $*$ qui est un difféomorphisme
involutif entre, au choix:
\begin{itemize}
\item Les polyèdres compacts hyperboliques et les polyèdres
compacts de Sitter.
\item Les polyèdres compacts sphériques sur eux-mêmes.
\end{itemize}
On va se contenter de définir l'application dans le premier cas:
On se place dans l'espace de Minkowski de dimension $4$ et on
considère $\mathbf{H}^{3}$ comme la composante connexe supérieure
de la sphère unité de type temps, et $\mathbf{S}_{1}^{3}$ comme la
sphère unité de type espace (cf.\ref{cons}). Dans ce modèle,
\begin{itemize} \item tout point de $\mathbf{H}^{3}$ s'interprète
comme droite de type temps.
\item toute droite hyperbolique s'interprète comme plan de
$\mathbf{R}_{1}^{3}$ de signature $(+,-)$.
\item tout plan hyperbolique s'interprète comme hyperplan de $\mathbf{R}_{1}^{3}$ orthogonal à une droite de type espace.
\end{itemize}
Et réciproquement :
\begin{itemize} \item tout point de $\mathbf{S}^{3}_{1}$
s'interprète comme droite de $\mathbf{R}_{1}^{3}$ de type temps.

\item tout plan de Sitter de type espace (c'est à dire de métrique induite riemannienne), s'interprète comme hyperplan de $\mathbf{R}_{1}^{3}$ orthogonal à une droite de type temps.
\item toute droite de Sitter de type espace s'interprète comme
plan de $\mathbf{R}_{1}^{3}$ de signature (+,+).
\end{itemize}
L'orthogonalité définit ainsi une correspondance duale entre
respectivement les points, les plans et les droites de
$\mathbf{H}^{3}$ d'une part et les plans de type espace, les point
et les droites de type espace de l'espace de Sitter.

\begin{Prodef}
Soit $P$ un polyèdre compact hyperbolique. Par la correspondance
ci-dessus on définit le polyèdre $P^{*}$ dans l'espace de de
Sitter vérifiant :
\begin{itemize}
\item Les sommets de $P^{*}$ sont "orthogonaux" aux faces de $P$.
\item Les arêtes de $P^{*}$ sont "orthogonales" aux arêtes de $P$.
\item les faces de $P^{*}$ sont "orthogonales" aux sommets de $P$.
\end{itemize}
C'est un polyèdre compact de l'espace de de Sitter, vérifiant les
propriétés suivantes :
\begin{itemize}
\item La combinatoire de $P^{*}$ est duale de celle de $P$.
\item La longueur des arêtes est égale aux angles dièdres duaux.
\item Réciproquement, les angles dièdres sont égaux aux longueurs
des arêtes duales.
\item Les angles entre arêtes adjacentes sont complémentaires aux
angles des arêtes adjacentes duales.

\end{itemize}
Enfin, la transformation $*$ commute aux isométries, ce qui
signifie que si $P$ est défini à isométrie prés, $P^{*}$ l'est
aussi.
\end{Prodef}
\paragraph {Remarques.}On peut construire la métrique induite sur $P^{*}$ en identifiant les arêtes communes des
polygônes sphériques duaux des sommets de $P$.\\
La somme des angles aux sommets de $P^{*}$ est supérieure à
$2\pi$, formant une singularité conique à courbure singulière négative.\\
Appliquer la formule de Gauss, et remarquer que la courbure
négative contenue dans les faces de $P$ se réalise comme
singularité conique négative dans les sommets de $P^{*}$, et
réciproquement, la courbure positive contenue dans les
singularités coniques des sommets de $P$ se réalise comme faces de
courbure constante positive.

\subsection{Le résultat de Rivin et Hodgson }
Rappel sur la topologie utilisée sur les espaces de polyèdres.
Lorsqu'on considère un espace de polyèdres, on peut définir une
topologie dessus de plusieures manières équivalentes :
\begin{itemize}
\item On considère la topologie Hausdorff définie sur les sous-ensembles de
$\mathbf{H}^{3}$.
\item On regarde les polyèdres dans l'espace projectif (même
topologie que $\mathbf{H}^{3}$ par continuité de la transformation
projective), et on paramètre les polyèdres par les sommets ou les
faces. L'espace de polyèdres est alors un sous-ensemble de
$\mathbf{R}^{3n}$, muni de la topologie induite, où $n$ est le
nombre de sommets ou de faces.
\end{itemize} On définit les espaces dans lesquels on va travailler :
\begin{Def}
On appellera $\mathcal{P}_{n}^{c}$ l'ensemble des polyèdres
compacts avec n faces numérotées, définis à une isométrie
hyperbolique prés. On notera de même $\mathcal{P}_{n}^{si}$
lorsque les polyèdres considérés sont semi-idéaux, et
$\mathcal{P}_{n}^{i}$ lorsque qu'ils sont idéaux .\end{Def} On
énonce alors :
\begin{Pro}
$\mathcal{P}_{n}^{c}$ est une variété différentielle de dimension
$3n-6$.Les espaces $\mathcal{P}_{n}^{si}$ et $\mathcal{P}_{n}^{i}$
sont inclus dans l'adhérence de $\mathcal{P}_{n}^{c}$.
\end{Pro}
\paragraph{Démonstration:}La condition de compacité et la condition d'existence de n faces sur
l'intersection des n demi-espaces est une condition ouverte sur
l'espace des n plans orientés numérotés formant les faces du
polyèdre. L'espace de polyèdres compacts à n faces numérotées est
donc un ouvert de $\mathbf{R}^{3n}$. Les isométries hyperboliques
forment un groupe de Lie de dimension 6 qui agit librement et
proprement sur cet ouvert. Le quotient est donc bien une variété
différentielle de la dimension souhaitée.\\On peut voir les
polyèdres semi-idéaux comme limite dans le modèle projectif de
polyèdres compacts de même combinatoire et dont les sommet tendent
vers la frontière de $\mathbf{H}^{3}$. $\diamond$

On définit maintenant les espaces métriques qui vont caractériser
les polyèdres compacts et semi-idéaux.
\begin{Def} On notera $\mathcal{S}_{n}$ l'espace des espaces
métriques $h$, tel que h soit homéomorphe à la sphère
$\mathbf{S}^{2}$, de courbure constante $1$ partout à l'exception
d'un nombre fini de points où il présente une singularité conique.
Les singularité sont numérotées, et les espaces métriques sont
définis aux isométries préservant les singularités prés.
\end{Def}
\paragraph{Remarque.}Si $P$ est un polyèdre hyperbolique compact, $P^{*}$ vu comme un espace
muni de la métrique induite est un élément de $\mathcal{S}_{n}$.
\begin{Def}
On notera $\mathcal{M}_{n}$ le sous-ensemble de $\mathcal{S}_{n}$
des métriques telles que :
\begin{itemize}
\item La courbure aux points singuliers est strictement négative.
\item Toute géodésique fermée a une longueur strictement supérieure
à $2\pi$. Ces métriques seront dites \emph{admissibles}.
\end{itemize}
On notera $\mathcal{M}^{si}_{n}$ le sous-ensemble de
$\mathcal{S}_{n}$ des métriques telles que:
\begin{itemize}
\item La courbure aux points singuliers est strictement négative.
\item toute géodésique fermée a une longueur
supérieure ou égale à $2\pi$.
\item le nombre de géodésiques fermées de longueur $2\pi$ est fini
et chacune d'entre elles sépare l'espace en deux composantes
connexes dont l'une au moins est isométrique à un hémisphère.
\end{itemize}
Ces métriques seront dites \emph{admissibles de volume fini}.
\end{Def}
Si $h_{1}$ et $h_{2}$ sont deux éléments de $\mathcal{M}_{n}$, on
considère les homéomorphismes $h_{1}\leftrightarrow h_{2}$
envoyant la $i$-ème singularité de $h_{1}$ sur la $i$-ème
singularité de $h_{2}$. Il existe une isométrie parmi ces
homéomorphismes, si et seulement si $h_{1}=h_{2}$. Sinon, on peut
considérer la borne inférieure des taux d'écart à une isométrie de
ces homéomorphismes. Ceci définit une métrique dite de Lipschitz
sur $\mathcal{S}_{n}$ (cf. \cite{comp} section 6). \\
Nous énonçons sans démonstration (cf. \cite{comp} et \cite{idea})
:
\begin{Pro}
$\mathcal{M}_{n}$ est une variété topologique de dimension 3n-6.
$\mathcal{M}^{si}_{n}$ est contenu dans l'adhérence de
$\mathcal{M}_{n}$.
\end{Pro}
A tout polyèdre compact de $\mathcal{P}_{n}$ on associe son dual
$P^{*}$ muni de sa métrique induite. En tant qu'objet intrinsèque,
c'est un élément de $\mathcal{S}_{n}$. On note $\Phi$
l'application ainsi définie. On est en mesure de citer le résultat
de \cite{comp} et \cite{idea}:
\begin{The}[Rivin et Hodgson]\label{Rivi}
\begin{enumerate}
\item L'application $\Phi$ ainsi définie est un homéomorphisme de
$\mathcal{P}_{n}$ sur $\mathcal{M}_{n}$.
\item Le prolongement par continuité de $\Phi$ à $\tilde{\Phi}:\mathcal{P}^{si}_{n}\rightarrow \mathcal{S}_{n}$ a
exactement pour image $\mathcal{M}^{si}_{n} $.
\end{enumerate}
Autrement dit:
\begin{enumerate}
\item Les polyèdres compacts sont en bijection avec les métriques
duales admissibles.(Existence et unicité).

\item les polyèdres semi-idéaux ont exactement pour métriques
duales les métriques admissibles de volume fini.(Existence).
\end{enumerate}
\end{The}
Les sections 4 et 5 de ce mémoire consiste à démontrer l'unicité
du second énoncé du théorème de Rivin, et d'en déduire une
caractérisation des polyèdres hyperidéaux par leurs angles
dièdres.
\subsection{Unicité des polyèdres semi-idéaux de
métriques duale donnée} On rappelle la version infinitésimale de
la rigidité des polyèdres euclidiens due à Aleksandrov
\cite{alex}:
\begin{The}[Rigidité des polyèdres euclidiens]
Toute déformation isométrique d'un polyèdre euclidien a pour
origine une isométrie infinitésimale.
\end{The}
\paragraph{Remarque.}Contrairement à la version du théorème de
rigidité cité dans l'introduction, la combinatoire n'est pas
fixée. Dans la démonstration, on utilise un argument
\textbf{abstrait}, pour montrer que la combinatoire est rigidifiée
si la métrique induite est fixe. Ce point est important, car il
est le point clef qui empêche de déduire la rigidité en fonction
des angles dièdres.

Afin de démontrer la rigidité des polyèdres semi-idéaux, on va
montrer que $\tilde{\phi}$ est localement injective, puis que
c'est un homéomorphisme local, et enfin un homéomorphisme global.
On est en mesure désormais d'énoncer le résultat :
\begin{Pro}\label{diff} L'application
$\tilde{\phi}:\mathcal{P}_{n}^{si} \rightarrow
\mathcal{M}^{si}_{n}$ est localement injective. Autrement dit on a
la rigidité locale des polyèdres semi-idéaux pour leur métriques
duale.
\end{Pro}

On a besoin du lemme suivant :
\begin{Lem}\label{conv}
Soit $P$ un polyèdre semi-idéal convexe. Dans le modèle de
Minkowski (cf. \ref{cons}), il existe un vecteur normal de type
espace, définissant un demi-espace ouvert de $\mathbf{R}_{1}^{4}$,
où $P^{*}$ est tout entier contenu.
\end{Lem}
\paragraph{Démonstration.}Soit p un point d'une face $F$ ouverte de
$P$, on lui associe sa droite $D$ de $\mathbf{R}^{4}_{1}$, que
l'on peut orienter par convexité. Je dis que le demi-espace ouvert
défini par $\perp D$ contient tout entier $P^{*}$ sauf le sommet
$F^{*}$ qui appartient à sa frontière. En effet par convexité,
$P^{*}$ est contenu dans le demi-espace fermé. Ensuite, le point
$p$ étant dans l'intérieur de $F$, $\perp D$ ne peut couper les
arêtes et les faces adjacentes à
$F^{*}$ qu'au niveau du sommet $F^{*}$.\\
Enfin, comme $F^{*}$ est distinct de l'origine de
$\mathbf{R}_{1}^{4}$, on peut perturber l'hyperplan $\perp D$ pour
obtenir le résultat souhaité.$\diamond$

\paragraph{Démonstration de \ref{diff}.} On va montrer que toute suite de
polyèdres de $\mathcal{P}^{si}_{n}$ convergeant vers $P$ et de
même métrique duale $h$ est stationnaire, ce qui impliquera
l'injectivité locale. Pour cela, on va raisonner par l'absurde et
supposer l'existence d'une suite ($P_{k}$) de polyèdres de
$\mathcal{P}^{si}_{n}$ convergeant vers $P$, tous les éléments
étant distincts de $P$. La dualité commutant aux isométries, on
obtient une suite $P_{k}^{*}$ de polyèdres (toujours définis à
isométrie prés) convergeant vers $P^{*}$, tous les éléments étant
distincts de $P^{*}$. On peut voir cette suite d'éléments comme
une suite de la variété différentielle de dimension $3n-6$
constituée de l'ensemble de n plans hyperboliques quotientés par
le groupe des isométries de $\mathbf{H}^{3}$. Un voisinage de $P$
dans cet espace est difféomorphe à un ouvert de
$\mathbf{R}^{3n-6}$. Par compacité de la sphère unité de
$\mathbf{R}^{3n-6}$, on peut extraire une sous-suite
$P^{*}_{s(k)}$ différentiable de polyèdres, c'est à dire
convergeant dans une direction donnée (les polyèdres étant vus
dans une carte comme vecteurs de $\mathbf{R}^{3n-6}$, la suite
$\frac{P_{s(k)}-P}{\parallel P_{s(k)}-P\parallel}$ converge vers
un champ de vecteur non nul sur $P$, noté $X$). La suite des
$P_{k}^{*}$ ayant engendré la déformation étant isométriques, le
champ $X$ est continu sur $P^{*}$ et est dérivable au sens de Lie
de dérivée de Lie nulle. On peut donc obtenir une "déformation
infinitésimale isométrique $X$ de $P^{*}$ non nulle au premier
ordre".

Maintenant, d'après \ref{conv}, on peut prendre un représentant
$P'^{*}$ de $P^{*}$ est dans $\Omega_{0}$, le demi-espace de
Sitter situé au dessus d'un certain hyperplan $H$ de type espace.
On choisit un champ de déformation isométrique (de dérivée de Lie
nulle) $Y$ de $P'^{*}$, qui se réalise à isométrie prés comme la
déformation infinitésimale $X$ de $P^{*}$. $X$ est non nul, ce qui
empêche $Y$ d'être la restriction d'une isométrie infinitésimale.
On applique à $P'^{*}$, l'application de Pogorelov infinitésimale
centrée sur la normale à l'hyperplan $H$. Le caractère polyèdral
et convexe étant un invariant projectif, on obtient, d'après le
théorème \ref{Pogi}, un polyèdre euclidien convexe muni d'une
déformation infinitésimale isométrique non nulle qui n'est pas la
restriction d'une isométrie  infinitésimale. Ceci contredit le
théorème de Cauchy.
Absurde.$\diamond$\\
Nous obtenons maintenant :
\begin{Pro} L'application
$\tilde{\phi}:\mathcal{P}_{n}^{si} \rightarrow
\mathcal{M}^{si}_{n}$ est un homéomorphisme local.
\end{Pro}
\paragraph{Démonstration.}$\mathcal{P}_{n}^{si}$ est un espace
topologique localement compact car c'est un sous-ensemble d'une
variété topologique, défini par un nombre fini de condition
ouverte et fermée. Or toute application localement injective
continue définie sur un espace localement compact est un
homéomorphisme local.$\diamond$

\begin{The}\label{The1} L'application
$\tilde{\phi}:\mathcal{P}_{n}^{si} \rightarrow
\mathcal{M}^{si}_{n}$ est un homéomorphisme global. On a donc une
relation bi-univoque entre les polyèdres semi-idéaux et les
métriques duales admissibles de volume fini.
\end{The}
\paragraph{Démonstration.}Il ne reste plus qu'à montrer
l'injectivité de: $\tilde{\phi}$ dans le cas où il y a au moins un
sommet idéal.\\
Un simple argument topologique va nous permettre de conclure.
Raisonnons par l'absurde et supposons que $h \in
\mathcal{M}^{si}_{n}$ (strictement semi-idéal), ait deux
antécédents distincts $P_{1}$ et $P_{2}$. Par séparation, il
existe deux voisinages $V_{1}$ de $P_{1}$ et $V_{2}$ de $P_{2}$
dans $\mathcal{P}_{n}^{si}$ disjoints. $\tilde{\phi}$ est un
homéomorphisme local, donc il existe un voisinage V de h dans
$\mathcal{M}^{si}_{n}$ et deux voisinages $W_{1}\subset V_{1}$ de
$P_{1}$ et $W_{2}\subset V_{2}$ de $P_{2}$ dans
$\mathcal{P}_{n}^{si}$ tel que $\tilde{\phi}$ soit un
homéomorphisme entre $W_{1}$ et $V$ et entre $W_{2}$ et $V$.Comme
$P_{1}$ et $P_{2}$ appartiennent à l'adhérence de
$\mathcal{P}_{n}^{c}$ et que $V$ est ouvert, $\tilde{\phi}$ ne
peut pas être injective sur $\mathcal{P}_{n}^{c}$,
absurde.$\diamond$
\section{Une caractérisation des polyèdres hyperboliques
hyperidéaux} Cette section est consacrée à la démonstration du
résultat de \cite{hype} comme corollaire de la caractérisation des
polyèdres semi-idéaux.
\subsection{Enoncé et définitions}
Pour définir les polyèdres hyperidéaux, on va considérer le modèle
projectif de l'espace hyperbolique $\mathbf{H}^{3}$. Dans ce
modèle $\mathbf{H}^{3}$ est représenté par le sous-ensemble
$\mathbf{R}^{3} \subset \mathbf{RP}^{3}$ constitué des droites
vectorielles de $\mathbf{R}^4_{1}$ de type temps.\\
D'autre part on définit les polyèdres de $\mathbf{RP}^{3}$:
\begin{Def}[Polyèdres projectifs.]
$P$ est un polyèdre projectif s'il existe une représentation
affine $\mathbf{R}^{3} \subset \mathbf{RP}^{3}$ de l'espace
projectif où $P$ est un polyèdre compact.
\end{Def}
\paragraph{Remarque.}Il existe des représentations affines de l'espace projectif où
$P$ est partiellement à l'infini.
\begin{Def}
On appelle polyèdre hyperidéal de $\mathbf{H}^{3}$ la donnée de
$n$ plans hyperboliques orientés formant un polyèdre projectif
comme plans de $\mathbf{RP}^{3}$ et dont les arêtes intersectent
$\mathbf{H}^{3} \subset \mathbf{RP}^{3}$ mais pas les sommets. Si
aucun des sommets n'appartient à la sphère à l'infini de
$\mathbf{H}^3$, le polyèdre est dit hyperidéal strict. Si au
contraire tous ses sommets s'y trouvent, il est dit idéal.
\end{Def}
On note l'application différentiable $\tilde{\Theta}$ qui à n
plans orientés hyperboliques associe leurs angles dièdres. Elle
est invariante par isométrie hyperbolique.
\begin{Def}On note $\Theta:\mathcal{P}_{\Gamma}\rightarrow ]0,\pi[^{E}$
l'application qui à tout polyèdre hyperidéal de combinatoire
$\Gamma$ associe la donnée de ses $E$ angles dièdres.
\end{Def}
Nous sommes en mesure maintenant d'énoncer le résultat à
démontrer. C'est le résultat principal de la prépublication de Bao
et Bonahon \cite{hype}, où une preuve directe est proposée.
\begin{Def}\label{defKgamma}
$K_{\Gamma}$ est le
sous-ensemble de $]0,\pi[^{E}$ défini de la manière suivante:\\
A chaque arête $e_{i}$ de $\Gamma^{*}$ on associe un poids
$\theta_{e_{i}}$ (représentant l'angle dièdre extérieur ie si
l'angle dièdre vaut $\alpha_{e} \in ]0,\pi[$ , l'angle dièdre
extérieur vaut $\theta_{e}= \pi-\alpha_{e} \in ]0, \pi[$) et
vérifiant les conditions suivantes :
\begin{itemize}
\item $\theta_{e_{i}} \in ]0,\pi[$
\item \emph{C1}- $\sum_{i=1}^{n}\theta_{e_{i}}\geq 2\pi$ pour toute courbe
fermée $\gamma$ plongée dans le graphe $\Gamma^{*}$ et passant par
$e_{1},...,e_{n}$; l'égalité ayant lieu très exactement quand
$\gamma$ est le bord d'une face du graphe dual $\Gamma^{*}$
correspondant à un sommet idéal.
\item \emph{C2}- $\sum_{i=1}^{n}\theta_{e_{i}} > \pi$ pour tout chemin $\gamma$ plongé
dans $\Gamma^{*}$ et joignant deux sommets de $\Gamma^{*}$ ayant
une même face adjacente, mais tel que $\gamma$ ne soit pas tout
entier contenu dans le bord de cette face.
\end{itemize}
\end{Def}
\begin{The}[Caractérisation des polyèdres hyperidéaux.] \label{caracthyp}L'application $\Theta$ est un homéomorphisme de
$\mathcal{P}_{\Gamma}$ vers $K_{\Gamma}$.

\end{The}
\subsection{L'application troncature}
On va définir un plongement différentiable de
$\mathcal{P}_{\Gamma}$ vers les polyèdres compacts, qui va nous
permettre de déduire le théorème \ref{caracthyp} du théorème
\ref{Rivi}.\\
On reprend notre modèle habituel où $\mathbf{RP}^{3}$ est
l'ensemble des droites vectorielles de l'espace de Minkowski
$\mathbf{R}_{1}^{4}$. A tout point de $\mathbf{RP}^{3}$ hors de
$\mathbf{H}^{3}$, ie à toute droite vectorielle de type espace, on
peut associer l'hyperplan orthogonal, qui s'identifie à un plan
hyperbolique que l'on note $\perp x$. Dans une représentation
affine de l'espace projectif, on obtient $\perp x$ en considérant
les droites affines passant par $x$ et tangentes à la sphère à
l'infini de $\mathbf{H}^{3}$, elles dessinent le bord à l'infini
de $\perp x$.\\
On peut caractériser $\perp x$ comme étant l'unique plan
hyperbolique perpendiculaire à toutes les droites hyperboliques
passant par $x$. En effet, si $d$ est une droite hyperbolique
passant par x dans l'espace projectif, cela signifie très
exactement qu'en tant que plan vectoriel, elle contient la droite
vectorielle $ x$. Or $x$ est l'orthogonal à $\perp x$. Le plan
vectoriel $d$ et l'hyperplan $\perp x$ sont donc perpendiculaires
dans l'espace de Minkowski; ce qui est équivalent à leur
perpendicularité comme objets hyperboliques. En effet,
$\mathbf{H}^{3}$ étant une surface de niveau du produit scalaire
dans
 l'espace de Minkowski, deux droites hyperboliques s'intersectant en un point sont orthogonales si et seulement si les plans correspondant dans $\mathbf{R}^{3}_{1}$ sont perpendiculaires. \\
 On peut maintenant
considérer un polyèdre hyperidéal $P$. A chacun de ses sommets
hyperidéaux $x$, on associe le demi-espace hyperbolique $H_{x}$.
De frontière $\perp x$, et orienté de sorte que x ne soit pas dans
le demi-espace affine associé. On énoncer la propriété :
\begin{Prodef}\label{trun}
Soit $P$ un polyèdre hyperidéal ayant $n$ faces et $\{x_{i}\}_{i
\in [1,k] }$ pour sommets hyperidéaux. On tronque $P$ en prenant
l'intersection avec les $H_{x_{i}}$ fournit un plongement
 $trun:\mathcal{P}^{h}_{n} \rightarrow \mathcal{P}^{si}_{n+k}$ dans l'espace des polyèdres compacts ayant $n+k$ faces.\\ Tous les plans hyperboliques $\perp x_{i}$ sont deux à deux
disjoints, ce qui détermine la combinatoire de
$trun (P)$.\\
Enfin, chaque $\perp x_{i}$ est caractérisé comme étant l'unique
plan hyperbolique perpendiculaire à toutes les arêtes adjacentes à
$x_{i}$.
\end{Prodef}
On a besoin du lemme suivant :
\begin{Lem}\label{coupeh3}Soit x et x' deux sommets hyperidéaux du polyèdre $P$.
La droite projective passant par ces deux points coupe
$\mathbf{H}^{3}$.
\end{Lem}
\paragraph{Démonstration.}On considère une représentation affine
du modèle projectif où $P$ est compacte. On considère la
projection radiale $\phi:\mathbf{R}^{3}-\{v \} \rightarrow
\mathbf{S}^{2}$ sur une sphère de centre x. L'image par $\phi$ de
l'espace hyperbolique est un disque sphérique ouvert. Par
convexité l'image par $\phi$ du polyèdre est contenu dans ce
disque. La droite $(xx')$ coupe donc ce disque et par  conséquent
$\mathbf{H}^{3}$.$\diamond$

\paragraph{Démonstration de la proposition-définition \ref{trun}.}Montrons tout d'abord que les $\perp
x_{i}$ sont deux à deux disjoints. Soit x et x' deux sommets
hyperidéaux. L'intersection des deux plans hyperboliques $\perp x$
et $\perp x'$ est un plan vectoriel de l'espace de Minkowski noté
$\Pi$ . L'orthogonal de $\Pi$ est très exactement dans le modèle
projectif, la droite passant par x et x'. D'après le lemme
\ref{coupeh3}, c'est un plan vectoriel de l'espace de Minkowski de
signature (+,-). Donc $\Pi$ est de type espace et ne coupe pas
l'espace hyperbolique.CQFD.\\

 $\perp x$ est  l'unique perpendiculaire commune aux arêtes
adjacentes: l'unicité découle du fait qu'il y a au moins
deux arêtes.\\

Il reste à montrer que $trun$ est un plongement. Pour cela on peut
maintenant facilement caractériser l'image : n étant fixé, c'est
l'ensemble des polyèdres compacts à $n+k$ faces vérifiant la
condition suivante:
\begin{itemize}
\item Il existe $k$ faces non adjacentes contenant tous les
sommets non idéaux, et auxquelles les arêtes adjacentes sont
perpendiculaires.
\end{itemize}
On reconstitue alors de manière différentiable le polyèdre initial
en prolongeant les arêtes de chaque face qui se coupent dans
l'espace projectif au point orthogonal de la face.$\diamond$

\subsection{Caractérisation des polyèdres hyperidéaux en terme de leur "métrique duale"}
Dans toute la suite, on se donne une combinatoire de polyèdre
$\Gamma$ avec $F$ faces numérotées, $A$ arêtes numérotées et $S$
sommets numérotés. Les sommets hyperidéaux $h$ seront distingués
des sommets idéaux par respectivement les index $h$ et $i$. Ainsi
$S_{h}$ et $S_{i}$ sont respectivement le nombre de sommets
hyperidéaux et idéaux. On considèrera $\mathcal{P}_{\Gamma}^{h}$,
l'ensemble des polyèdres hyperidéaux de combinatoire $\Gamma$
fixée.  On a un homéomorphisme par composition des applications
troncature et $\tilde{\phi}$, qui à un polyèdre semi-idéal associe
sa métrique duale:
$$trun \circ \tilde{\phi}: \mathcal{P}_{\Gamma}^{h} \rightarrow
\mathcal{M}_{S_{h}+F}^{si}$$ On notera
$\mathcal{M}^{h}_{\Gamma}=trun \circ
\tilde{\phi}(\mathcal{P}_{\Gamma}^{h})$ l'image de l'application.
Le but de cette sous-section est de caractériser cet ensemble.\\
\begin{Def}["néga-et posi-hémisphère"]On appelle
\emph{"néga-hémisphère"}, (respectivement
\emph{"posi-hémisphère"}), la variété à bord, riemannienne de
courbure constante $+1$, avec une singularité conique de courbure
négative (respectivement positive), homéomorphe à un disque fermé,
et tel que le bord soit exactement l'ensemble des points situés à
distance $\frac{\pi}{2}$ de la singularité.

\end{Def}
On énonce quelques propriétés de ces variétés :
\begin{Pro}
Les néga- et posi- hémisphère sont déterminés uniquement par la
valeur de leur singularité qui va de 0 à l'infini.\\
On obtient le néga- ou posi- hémisphère de singularité conique
$\alpha$ de la manière suivante:
\begin{itemize}
\item On prend $n$ triangles sphériques doublement rectangles d'angle
au sommet $\alpha_{i}$ . Il a donc deux côtés de longueur
$\frac{\pi}{2}$, et un côté de longueur $\alpha_{i}$. On demande
que la somme des $\alpha_{i}$ vaille $\alpha$.
\item On identifie isométriquement les côtés de longueur
$\frac{\pi}{2}$des triangles $i$ et $i+1$ (modulo n).
\end{itemize}
La longueur du bord vaut $ \alpha$. \end{Pro}
\begin{Pro}\label{hemipi} Toute géodésique complète d'un néga- ou posi- hémisphère est un
segment immergé de longueur $\pi$ dont les extrémités sont dans le
bord de la variété. Ces géodésiques ne se recoupent jamais si la
variété est un néga-hémisphère.
\end{Pro}

\paragraph{Démonstration.}Les étapes :\\
Caractérisation par la valeur de la singularité : L'"espace
tangent" au niveau de la singularité est un cône circulaire
d'angle $\alpha$. On peut définir l'équivalent d'une application
exponentielle bijective en prenant les vecteurs du cône de norme
inférieure ou égale à $\frac{\pi}{2}$. Ce geste permet de
construire par composition une isométrie entre posi- ou néga-
hémisphères de même singularité.\\
La construction par des triangles vérifie bien la définition d'un
posi- ou néga- hémisphère.\\
Longueur des géodésiques.  Les géodésiques passant par une
singularité sont toujours l'union de deux arcs de longueur
$\frac{\pi}{2}$ joignant la singularité et le bord. Dans le cas où
la géodésique ne passe pas par la singularité, on va montrer que
sa longueur est $\pi$ d'abord pour les négahémisphères. Soit une
telle géodésique. Découpons le négahémisphère en enlevant un
triangle bi-rectangle de sommet principal le pôle et de côté
opposé situé sur le bord, de telle sorte à obtenir un hémisphère
simple. On s'arrange pour conserver un morceau de géodésique non
inclus dans le bord. Ce dernier appartient sur l'hémisphère simple
à un grand cercle qui n'est pas l'équateur . La géodésique
initiale est donc le prolongement de ce morceau, ie un segment
immergé de longueur $\pi$, dont les extrémité sont sur le bord.
Pour le posi-hémisphère, on le découpe du bord à la singularité
puis on le recolle isométriquement suffisamment de fois avec
lui-même afin d'obtenir un néga-hémisphère. On applique alors le
résultat précédent, et on en déduit que la géodésique est
un segment immergé de longueur $\pi$ qui "s'enroule" autour du posihémisphère. $\diamond$\\
On va maintenant pouvoir caractériser $\mathcal{M}_{\Gamma}$. Pour
cela on va définir un certain ensemble d'espaces métriques :
\begin{Def}\label{defemi}On se donne la combinatoire $\Gamma$ d'un
polyèdre hyperidéal. A chaque arête on associe un poids
$\theta_{e_{i}}$ correspondant à l'angle dièdre extérieur, de
telle sorte que la somme autour d'un sommet hyperidéal soit
supérieure à $2\pi$, et la somme autour d'un sommet idéal soit
égale à $2\pi$. On construit un espace métrique de la manière
suivante :
\begin{itemize}
\item On prend autant de négahémisphères que de sommets
hyperidéaux, la valeur de la singularité valant la somme des
angles dièdres extérieurs autour du sommet. On prend autant
d'hémisphères droits que de sommets idéaux.
\item On place sur le bord de chaque (néga-)hémisphère, dans
l'ordre, autant de singularités qu'il y a de faces adjacentes au
sommet en question dans $\Gamma$. La distance entre chaque
singularité valant l'angle dièdre extérieur correspondant.
\item On identifie isométriquement les segments ainsi construits
ayant la même arête pour origine dans $\Gamma$.
\end{itemize}
On demande que toutes les géodésiques fermées de cet espace soient
strictement supérieures à $2\pi$, sauf uniquement s'il s'agit du
bord d'une hémisphère simple.\\
Les singularité situées sur le bord des hémisphères sont des
singularités de courbure strictement négative multiples de $\pi$.\\
 L'ensemble de métriques avec $S_{h}+F$ singularités numérotées
ainsi obtenus, sous-ensemble de $\mathcal{M}^{si}_{S_{h}+F}$, sera
noté $\mathcal{Q}_{\Gamma}$. Les $S_{h}$ premières singularités,
pôles des néga-hémisphères, seront dite de type "h", les $F$
dernières seront dite de type "f".
\end{Def}
\paragraph{Remarque.}La notation (néga-)hémisphère signifie qu'il peut éventuellement s'agir d'un hémisphère simple.
 On énonce maintenant la caractérisation attendue :
\begin{The}\label{caracmet}
On a l'égalité $\mathcal{Q}_{\Gamma}=\mathcal{M}_{\Gamma}$.\\
Un espace métrique $h$ de cet ensemble avec ses $S_{h}+F$
singularités numérotées étant donnée, on sait retrouver
directement les géodésiques de $h$ formant la combinatoire duale
du polyèdre tronqué, image de $h$ par $\tilde{\phi}^{-1}$.
\end{The}
\paragraph{Démonstration.}Montrons tout d'abord la première inclusion:
$\mathcal{M}_{\Gamma}\subset \mathcal{Q}_{\Gamma}$. En effet, soit
$P$ un polyèdre hyperidéal tronqué. Les "faces" du polyèdre dual
correspondant aux sommets idéaux sont isométriques à des
hémisphères simples. Considérons maintenant les sommets d'une même
face tronquée. Par perpendicularité de cette face par rapport aux
faces adjacentes, les faces du polyèdre dual correspondant aux
sommets sont des triangles sphériques bi-rectangles, collés les
uns aux autres pour former un néga-hémisphère de singularité
conique de valeur l'aire de la face tronquée plus $2\pi$. On
obtient ainsi toutes les faces du polyèdre dual agencées selon la
combinatoire duale $\Gamma^{*}$. Il reste à montrer que les
longueurs des géodésiques fermées sont strictement supérieures à
$2\pi$, sauf uniquement si elles s'identifient au bord d'un
hémisphère simple. C'est une conséquence
directe du théorème \ref{The1}.\\
Soit $h \in \mathcal{M}_{\Gamma}$. On note $\Lambda$ l'ensemble
des géodésiques formant la combinatoire du polyèdre original
Montrons que l'on peut récupérer directement cette combinatoire
par le processus $1$ suivant : considérer une singularité de type
"f" situé au point $A$, et sa singularité la plus proche située au
point $B$. Par construction de $\mathcal{M}_{\Gamma}$, il existe
une géodésique joignant $A$ et $B$, incluse dans $\Lambda$ et
strictement inférieure à $\pi$. Choisissons une géodésique
$\gamma$ joignant $A$ à $B$, et de longueur inférieure à $\pi$.
Comme $h$ est constitué de recollement de (néga-)hémisphères,
d'après la proposition \ref{hemipi}, cette géodésique appartient à
$\Lambda$. Ensuite, au niveau de la singularité $B$, on trace
toutes les géodésiques partant de $B$ et formant un angle multiple
de $\pi$ avec la géodésique $AB$. Elles sont dans $\Lambda$. On
réitère ce procédé et par connexité de $\Gamma^{*}$, au bout d'un
nombre fini d'étapes, on peut lire sur $h$ toute
la combinatoire $\Lambda$.\\
Montrons la dernière inclusion : $\mathcal{Q}_{\Gamma}\subset
\mathcal{M}_{\Gamma}$. Soit $g \in \mathcal{Q}_{\Gamma}$. D'après
le théorème \ref{The1}, résultat principal de la section
précédente, il existe un unique polyèdre semi-idéal $P$ ayant $g$
pour métrique duale.\\
Regardons sur $g$ la combinatoire $\Sigma^{*}$ induite par $P$. On
cherche à montrer que cette combinatoire est la même que la
combinatoire $\Gamma^{*}$ obtenue sur g par le
processus $1$ ( processus qui aboutit par construction de $g$).\\
Rappelons que les arêtes de $\Sigma^{*}$ sur $g$ correspondent à
des angles dièdres et sont donc strictement inférieures à $\pi$.
Soit $a$ une arête de $\Sigma^{*}$ non incluse dans $\Gamma^{*}$.
Si $a$ joint deux singularités de type "h", elle traverse deux
hémisphères et est supérieure ou égale à $\pi$. Si elle joint deux
singularités de type "f", elle fait un séjour dans un
(néga-)hémisphère et est supérieure ou égale à $\pi$ (par
propriété des géodésiques des (néga-)hémisphères). Si $a$ joint
deux singularités de type différents, elle passe par deux
hémisphères différents et pour la même raison est supérieure ou
égale à $\pi$. Donc $\Sigma^{*} \subset \Gamma^{*}$. Supposons que
l'inclusion soit stricte, on pourrait alors trouver deux arêtes de
$\Sigma^{*}$, adjacentes et bordant une même face avec un angle
supérieur à $\pi$. Ceci est interdit par
convexité du polyèdre $P$. Donc $\Sigma^{*} = \Gamma^{*}$. CQFD.\\
Le polyèdre $P$ possède donc la combinatoire souhaitée, les
sommets correspondant aux hémisphères droits de $g$ sont idéaux,
et les faces correspondant aux singularités de type "h" sont
perpendiculaires à leur faces adjacentes. On peut donc appliquer
$trun^{-1}$ et conclure.$\diamond$
\subsection{Caractérisation des polyèdres hyperidéaux en terme de leurs angles dièdres}
Le but de cette section est de remarquer que l'espace des
métriques duales des polyèdres hyperidéaux est suffisamment
dégénéré pour être décrit en terme de graphe muni de poids
correspondants aux angles dièdres.

 On rappelle la chose suivante :
\begin{Lem}
Soit $\mathcal{S}$ un sous-ensemble de $\mathcal{S}_{n}$, espaces
des métriques sphériques avec n singularités numérotées. On
demande que les éléments de $\mathcal{S}_{n}$ admettent tous une
même triangulation par $E$ géodésiques de longueurs strictement
inférieures à $\pi$.L'application $\psi:\mathcal{S} \rightarrow
]0,\pi[^{E}$, qui à un espace métrique associe la longueur des
géodésique de sa triangulation est un plongement.
\end{Lem}
\paragraph{Démonstration.}$\psi^{-1}$ est la restriction de
l'application continue et injective qui à un E-uplet de
$]0,\pi[^{E}$ associe l'espace obtenues en recollant les triangles
sphériques déterminés uniquement par la longueur de leurs
côtés.$\diamond$\\
\\
Soit $h \in \mathcal{Q}_{\Gamma}$, on a vu que l'on pouvait
retrouver la combinatoire $\Gamma^{*}$ sur $h$. On peut trianguler
les hémisphères correspondant aux sommets idéaux. On peut donc
définir un plongement $\psi:\mathcal{Q}_{\Gamma} \rightarrow ]0,
\pi[^{E}$. Inversement,(on rappelle que $K_{\Gamma}$ est le
sous-ensemble de $]O,\pi[^{E}$ défini au théorème
\ref{caracthyp}), on peut définir une application $\chi:K_{\Gamma}
\rightarrow \mathcal{S}_{S_{h}+F}$ définie comme dans la
construction de $\mathcal{Q}_{\Gamma}$ dans la définition
\ref{defemi} en enlevant la condition pour les géodésiques d'être
de longueur strictement supérieures à $2\pi$. Rappelons cette
construction : On recolle des hémisphères et des négahémisphères
avec les longueurs d'arêtes prescrites par $K_{\Gamma}$, selon la
combinatoire $\Gamma^{*}$. Les conditions sur $K_{\Gamma}$
assurent d'avoir des hémisphères au niveau des faces duales des
sommets de type "i", et des négahémisphères au niveau des faces
duales des sommets de type "h".\\
Sur l'ensemble où elle est définie l'application composée $\psi
\circ \chi$ est l'identité.
\begin{The}
$$\psi(\mathcal{Q}_{\Gamma})=K_{\Gamma}$$ ie
$$\chi(K_{\Gamma})=\mathcal{Q}_{\Gamma}$$
\end{The}
\paragraph{Démonstration.}

Montrons la première inclusion:\\
$\psi(\mathcal{Q}_{\Gamma}) \subset K_{\Gamma}$. Soit $h \in
\mathcal{Q}_{\Gamma}$ avec les géodésiques dessinant la
combinatoire. On veut donc montrer que la condition pour toutes
les géodésiques fermées de $h$ d'être strictement supérieures à
$2\pi$ implique que le graphe de géodésiques de combinatoire
$\Gamma^{*}$ (les géodésiques joignant une singularité de type "h"
aux singularités de type "f" ne sont pas prises en compte) vérifie
les conditions de la définition de
$K_{\Gamma}$ (théorème \ref{caracthyp}).\\
Première condition. Soit $\gamma$ une courbe fermée de
$\Gamma^{*}$. $\gamma$ est une géodésique de $h$ puisqu'à chaque
singularité, $\gamma$ est localement une géodésique (pas de "coin"
d'angle inférieur strictement à $\pi$). La condition sur les
géodésiques fermées de $h$ entraîne que $\gamma$ est de longueur
strictement plus grande que $2\pi$, sauf uniquement dans le cas où
$\gamma$ borde un hémisphère simple.
Ceci constitue exactement la première condition.\\
Seconde condition. Soit $\gamma$ un chemin sur $\Gamma^{*}$ dont
les extrémités sont dans le bord d'une même (néga-)hémisphère,
mais non contenu tout entier dans ce bord. On cherche à minorer la
longueur de $\gamma$ par $\pi$; quitte à enlever une partie de la
courbe, on peut supposer que seules les extrémités de $\gamma$
sont sur le bord de l'hémisphère. Maintenant on peut fermer la
géodésique en une géodésique fermée en joignant les deux extrémité
par un arc de géodésique appartenant au (néga-)hémisphère . Cet
arc sera contenu dans le bord si les extrémités sont proches, ou
au contraire passera par le pôle du (néga-)hémisphère si elles
sont éloignées. Il est de longueur inférieure ou égale à $\pi$. La
condition sur les géodésiques de
$h$ permet de minorer strictement la longueur de  $\gamma$ par $\pi$.CQFD.\\
Montrons maintenant l'autre inclusion: $\chi(K_{\Gamma}) \subset
\mathcal{Q}_{\Gamma}$. Soit donc $h \in \chi(K_{\Gamma})$. Soit
$\gamma$ une géodésique fermée de $h$. On cherche à montrer que sa
longueur est strictement supérieure à $2\pi$, connaissant
les deux conditions (C1) et (C2) imposées aux angles dièdres.\\
Ou bien $\gamma$ est contenu dans $\Gamma^{*}$, et le résultat est
immédiat, ou bien $\gamma$ possède un arc de longueur $\pi$ dans
un (néga-)hémisphère. Considérons maintenant le chemin $\gamma'$
constitué de $\gamma$ privé de cet arc. On remarque que $\gamma'$
ne peut être contenu dans le bord du (néga-)hémisphère (sinon
$\gamma$ ne serait pas géodésique au niveau des extrémités de
$\gamma'$). Ou bien $\gamma'$ est contenu dans $\Gamma^{*}$ et le
résultat s'en déduit d'après la condition (C2) de la définition
\ref{defKgamma} (cette condition n'est pas nécessaire si tous les
(néga-)hémisphères sont des hémisphères simples), ou bien $\gamma$
posséde deux arcs de géodésiques sur deux hémisphères différentes.
Maintenant ou bien $\gamma$ n'est pas la réunion de ces deux arcs
est sa longueur est supérieure strictement à $2\pi$, ou bien elle
l'est effectivement. Dans ce cas, les deux (néga-)hémisphères sont
en contact sur une longueur supérieure à $\pi$, ce qui est
impossible.$\diamond$
\paragraph{Remarque.} La condition (C2) est inutile dans le cas
idéal. La caractérisation de polyèdres idéaux en fonction de leurs
angles dièdres ne fait intervenir que la condition (C1).

\subsection{Conclusion}
On a donc existence et unicité des polyèdres hyperidéaux de
combinatoire $\Gamma$ dont les angles dièdres extérieurs vérifient
les deux
conditions suivante :\\
(\emph{C1})$\sum_{i=1}^{n}\theta_{e_{i}}\geq 2\pi$ pour toute
courbe fermée $\gamma$ plongée dans le graphe $\Gamma^{*}$ et
passant par $e_{1},...,e_{n}$; l'égalité ayant lieu très
exactement quand $\gamma$ est le bord d'une face du graphe
dual $\Gamma^{*}$ correspondant à un sommet idéal.\\
(\emph{C2}) $\sum_{i=1}^{n}\theta_{e_{i}} > \pi$ pour tout chemin
$\gamma$ plongé dans $\Gamma^{*}$ et joignant deux sommets de
$\Gamma^{*}$ ayant une même face adjacente, mais tel que $\gamma$
ne soit pas tout entier contenu dans le bord de cette face.

Cette caractérisation permet d'avoir un critère pour savoir si un
polyèdre de combinatoire fixée est inscriptible dans une sphère:\\
On a existence et unicité des polyèdres idéaux de combinatoire
$\Gamma$ dont les angles dièdres vérifient la condition suivante :\\
(\emph{C1})$\sum_{i=1}^{n}\theta_{e_{i}}\geq 2\pi$ pour toute
courbe fermée $\gamma$ plongée dans le graphe dual $\Gamma^{*}$ et
passant par $e_{1},...,e_{n}$; l'égalité ayant lieu très
exactement quand $\gamma$ est le bord d'une face de $\Gamma^{*}$.

\section{Le cas fuchsien}Dans cette section, on se propose de
généraliser tous les résultats précédents au cas fuchsien. Les
résultats de Rivin et Hodgson \cite{comp} et \cite{idea} ont été
généralisés au cas fuchsien par Schlenker dans \cite{fuchs}. Nous
nous proposons d'appliquer la même démarche que dans les sections
précédentes à ce cas.
\subsection{Définitions}
On sait que pour tout $g \geq 2$, il existe un sous-groupe discret
$\Gamma \subset Isom(\mathbf{H}^{2})$, tel que
$\mathbf{H}^{2}/\Gamma$ soit une surface compacte hyperbolique de
genre $g \geq 2$. On va étendre l'action de ce groupe à
$\mathbf{H}^{3}$ de la manière suivante : On considère un plan
hyperbolique $\mathbf{H}^{2} \subset \mathbf{H}^{3}$ et on écrit
$\mathbf{H}^{3}$ comme l'ensemble des couples $(t,x) \in
\mathbf{R} \times \mathbf{H}^{2}$ muni de la métrique
$(dt^{2}+ch^{2}(t)g_{\mathbf{H}^{2}})$. On fait alors agir $\gamma
\in  \Gamma$ par la formule : $\gamma(t,x)=(t,\gamma x)$.
\begin{Def}Le quotient $\mathbf{H}^{3}/\Gamma$ obtenu
(topologiquement le produit cartésien d'une surface compacte de
genre $g$ avec $\mathbf{R}$) est appelé variété fuchsienne.
\end{Def}
\begin{Def}Soit $S$ une surface fermée de genre $g \geq 2$.On note $\tilde{S}$
son revêtement universel (le plan) et $\Gamma$ sont groupe
fondamental. Un polyèdre compact convexe fuchsien est la donnée
d'un couple $(\phi,\rho)$ vérifiant les conditions suivantes :
\begin{itemize}
\item $\phi:\tilde{S} \rightarrow \mathbf{H}^{3}$ est un
plongement polyèdral convexe, i.e. la donnée d'une cellulation de
$\tilde{S}$ telle que $\phi$ envoie les arêtes de la cellulation
sur des segments compacts de géodésiques et les faces de la
cellulation sur des portions compactes de plans, formant un
ensemble polyèdral convexe.
\item $\rho:\Gamma \rightarrow Isom(P)$ est une
action de $\Gamma$ sur un plan $P \equiv \mathbf{H}^{2}$ de
$\mathbf{H}^{3}$, ie un morphisme de groupe dans l'ensemble des
isométrie de $\mathbf{H}^{3}$ laissant $P$ invariant.
\item $\phi$ est $\Gamma$-équivariant, ie :
$$\forall x \in \tilde{S},\forall \gamma \in \Gamma, \phi(\gamma
x)=\rho(\gamma)\phi(x)$$
\end{itemize}
En quotientant $\phi$ par $\Gamma$, on obtient un plongement
convexe polyèdral de $S$ dans la variété fuchsienne
$\mathbf{H}^{3}/\Gamma$.
\end{Def}
\paragraph{Remarques.}
\begin{itemize}

\item Le nombre de faces et de sommets de $\phi(\tilde{S})$ est infini.
\item Dans l'espace projectif, le bord de $\phi(\tilde{S})$ est
très exactement le bord de $P$.
\item La combinatoire $\Gamma$ d'un polyèdre fuchsien et sa
combinatoire duale sont des décompositions cellulaires d'une
surface de genre $g \geq 2$.
\end{itemize}
\begin{Def}On étend la définition précédente en considérant le cas
où le plongement polyèdral convexe $\phi$ est à valeurs dans
l'espace projectif $\mathbf{RP}^3$ muni du groupe des isométrie
hyperboliques. Un polyèdre fuchsien est alors dit respectivement
semi-idéal, idéal, hyperidéal, strictement hyperidéal,si dans
l'espace projectif toutes les arêtes du plongement polyèdral
$\phi$ coupent $\mathbf{H}^{3}$ et si respectivement, certains
sommets sont sur le bord de $\mathbf{H}^{3}$, tous les sommets
sont sur le bord de $\mathbf{H}^{3}$, les sommets sont soit sur le
bord soit au-delà de $\mathbf{H}^{3}$, tous les sommets sont
au-delà de $\mathbf{H}^{3}$.
\end{Def}
\subsection{Caractérisation des polyèdres fuchsiens semi-idéaux
en fonction de leur métrique duale} On se donne un polyèdre
fuchsien $(\phi,\rho)$ dans $\mathbf{H}^{3}$. Vu dans
$\mathbf{S}_{1}^{3}$, la représentation $\rho$, dont l'image est
un sous-groupe de transformations orthogonales de
$\mathbf{R}^{3}_{1}$, fixe un point (l'orthogonal du plan $P$). Le
plongement dual $\phi^{*}$ donne un plongement polyèdral convexe
de type espace, équivariant dans l'espace de Sitter, car les
images de la représentation $\rho$, transformations orthogonales,
commutent avec l'application $*$. La métrique duale du polyèdre
fuchsien est la métrique induite sur l'image de ce plongement
$\phi^{*}$. Comme souligné dans \cite{fuchs}, on a un résultat de
rigidité infinitésimale en fonction de la métrique duale.
\begin{Pro}\label{ishka}Soit $(\phi,\rho)$ un polyèdre fuchsien de
$\mathbf{H}^{3}$. Il n'y a pas de déformation infinitésimale
non-triviale (i.e. pas de la forme $\frac{d}{dt}(\epsilon_{t}
\phi,\epsilon_{t} \rho \epsilon_{t}^{-1})$ où $\epsilon_{t}$ est
une isométrie infinitésimale de $\mathbf{H}^{3}$), qui ne déforme
pas la métrique duale de l'image de $\phi$ au premier ordre.
\end{Pro}
\paragraph{Schéma de démonstration.}
Comme expliqué dans le lemme 4.8 de \cite{fuchs}, on utilise une
généralisation des idées de Cauchy sur la rigidité des polyèdres
aux surfaces équivariantes de genre supérieur ou égal à 2 dans
l'espace de Minkowski (cf. \cite{genren}). Puis on utilise la
transformation de Pogorelov dans l'espace de Sitter, qui
fonctionne à condition que $\rho$ laisse invariant un point (cf.
\cite{lorent}).$\diamond$\\

Maintenant, fixons le nombre $n$ de faces du polyèdres et le genre
$g$ de la surface $S$ considérée.
\begin{Def}
On appelle : $\mathcal{P}_{f}^{c}$ l'ensemble des polyèdres
fuchsiens de genre $g$, avec $n$ faces, définis à isométrie prés
(ie quotienté par la relation d'équivalence $(u\phi,u\rho
u^{-1})=(\phi,\rho)$ où $u$
est une isométrie directe laissant le plan $P$ invariant).\\
On note $\mathcal{P}_{f}^{si}$ l'ensemble des polyèdres fuchsiens
semi-idéaux (le plongement polyèdral $\phi$ présente des sommets à
l'infini), et $\mathcal{P}_{f}^{h}$ l'ensemble des polyèdres
fuchsiens hyperidéaux (le plongement polyèdral $\phi$ présente des
sommets hyperidéaux).
\end{Def}
On définit les espaces de métriques de manière analogue au cas de
genre $0$.
\begin{Def} On notera $\mathcal{M}_{f}^{c}$ l'ensemble des espaces
métriques $h$, définis à isométrie prés, tel que h soit
homéomorphe à la surface $S$ de genre $2$, de courbure constante
$1$, partout à l'exception d'un nombre fini égal à $n$ de points
où il présente une singularité conique strictement supérieure à
$2\pi$. Les singularités sont numérotées. On demande que les
géodésiques contractiles soient de longueur strictement
supérieures à $2\pi$.
\end{Def}
\begin{Def} On notera $\mathcal{M}_{f}^{i}$ l'espace des espaces
métriques $h$, définis à isométrie prés, tel que h soit
homéomorphe à la surface $S$ de genre $2$, de courbure constante
$1$ partout à l'exception d'un nombre fini égal à $n$ de points où
il présente une singularité conique d'angle strictement supérieure
à $2\pi$. Les singularités sont numérotées. On demande que les
géodésiques contractiles soient de longueur strictement supérieure
à $2\pi$, sauf un nombre fini qui sont alors le bord d'un
hémisphère.
\end{Def}
On énonce le théorème démontré dans la section 4 de \cite{fuchs}:
\begin{The}[Schlenker]\label{Schlenker}On a les deux propriétés
suivantes :
\begin{enumerate}
\item L'application $\phi$ qui à un polyèdre fuchsien associe sa métrique duale est un homéomorphisme de
$\mathcal{P}_{f}^{c}$ sur $\mathcal{M}^{c}_{f}$.
\item Le prolongement par continuité de $\phi$ à $\tilde{\phi}:\mathcal{P}_{f}^{si}\rightarrow \mathcal{S}_{f}$ a
exactement pour image $\mathcal{M}^{si}_{f} $.
\end{enumerate}
\end{The}
On est en mesure de récolter l'injectivité $\tilde{\phi}$, en
utilisant de manière analogue à la démonstration de la proposition
\ref{diff}, la rigidité infinitésimale. Ce point est absent de
\cite{fuchs}.
\begin{The}
L'application $\tilde{\phi}:\mathcal{P}^{si}_{f}\rightarrow
\mathcal{M}^{si}_{f} $ est localement injective.
\end{The}
\paragraph{Démonstration.} On suppose par l'absurde que l'application n'est pas
localement injective en un polyèdre fuchsien semi-idéal $P$. De la
même manière que dans la démonstration de \ref{diff}, on peut
extraire une suite différentiable de polyèdres semi-idéaux non
nulle au premier ordre, mais dont la métrique duale est constante.
Ceci est interdit par la proposition
\ref{ishka}.$\diamond$\\
On obtient de manière analogue au théorème \ref{The1}:
\begin{The}\label{The2} L'application
$\tilde{\phi}:\mathcal{P}_{f}^{si} \rightarrow
\mathcal{M}^{si}_{f}$ est un homéomorphisme global. On a donc une
relation bi-univoque entre les polyèdres fuchsien semi-idéaux et
leur métrique duale ie $\mathcal{M}^{si}_{f}$.
\end{The}
\paragraph{Démonstration.}
Elle est la même que dans le cas non-fuchsien (reprendre le
théorème \ref{The1}). En effet, l'espace des polyèdres fuchsiens
semi-idéal est localement compact (L'espace des polyèdres
fuchsiens $(\phi,\rho)$, à n faces dont les sommets sont
indifféremment compacts, idéaux, ou hyperidéaux, et défini à
isométrie laissant P invariant prés, est localement compact. Les
conditions de semi-idéalité sont une intersection finies d'ouverts
et de fermés dans cet espace). $\tilde{\phi}$ qui est localement
injective, est donc un homéomorphisme local. L'argument
topologique nécessaire pour conclure fonctionne car la situation
est identique au cas non fuchsien:
\begin{itemize}
\item $\mathcal{P}_{f}^{si}$ est dans l'adhérence de
$\mathcal{P}_{f}^{c}$ et $\mathcal{M}_{f}^{si}$ est dans
l'adhérence de $\mathcal{M}_{f}^{c}$.
\item $\mathcal{P}_{f}^{c}$ est homéomorphe à
$\mathcal{M}_{f}^{c}$.
\item $\mathcal{P}_{f}^{si}$ est localement homéomorphe à
$\mathcal{M}_{f}^{si}$.
\end{itemize}
$\diamond$\\
On peut maintenant se demander si on peut généraliser la
caractérisation des polyèdres hyperidéaux en fonction de leurs
angles dièdres au cas fuchsien.\\
\subsection{Caractérisation des polyèdres fuchsiens hyperidéaux en fonction de leur "métrique duale"}
On note $\mathcal{P}_{f,\Gamma}$ l'ensemble des polyèdres
fuchsiens hyperidéaux de combinatoire donnée $\Gamma$. Soit $P$ un
polyèdre fuchsien hyperidéal. Dans $\mathbf{H}^{3}$, on peut
tronquer les sommets strictement hyperidéaux du plongement
équivariant $\phi$ par des plans hyperboliques disjoints (cf.
proposition-définition \ref{trun}). L'application troncature est
donc une transformation purement géométrique qui "passe au
quotient". Si l'on cherche à caractériser maintenant l'image par
$\tilde{\phi}$ de ces polyèdres tronqués, on obtient le même
résultat que dans le théorème \ref{caracmet} moyennant quelques
modifications. On note $\mathcal{M}^{h}_{f,\Gamma}$ cette image.
On va définir tout d'abord l'espace de métriques candidat pour
être l'image de $\tilde{\phi}$:
\begin{Def}\label{defuchs} On se donne la combinatoire de genre $g$, $\Gamma$ d'un
polyèdre fuchsien hyperidéal. A chaque arête on associe un poids
$\theta_{e_{i}}$ correspondant à l'angle dièdre, de tel sorte que
la somme autour d'un sommet hyperidéal soit supérieure à $2\pi$,
et la somme autour d'un sommet idéal soit égale à $2\pi$. On
construit un espace métrique de la manière suivante :
\begin{itemize}
\item On prend autant de négahémisphères que de sommets
hyperidéaux, la valeur de la singularité valant la somme des
angles dièdres autour du sommet. On prend autant d'hémisphères
droits que de sommets idéaux.
\item On place sur le bord de chaque (néga-)hémisphère, dans
l'ordre, autant de points qu'il y a de faces adjacentes au sommet
en question dans $\Gamma$. La distance entre chaque points valant
l'angle dièdre correspondant.
\item On identifie isométriquement les segments ainsi construits
ayant la même arête pour origine dans $\Gamma$.
\end{itemize}
On demande que toutes les géodésiques contractiles fermées de cet
espace soient de longueur strictement supérieures à $2\pi$, sauf
uniquement s'il s'agit du
bord d'une hémisphère simple.\\
Les points situées sur le bord des hémisphères sont des
singularités de courbure strictement négative multiples de $\pi$.\\
L'ensemble de métriques avec $S_{h}+F$ singularités numérotées
ainsi obtenus, sous-ensemble de $\mathcal{M}^{si}_{f,S_{h}+F}$,
sera noté $\mathcal{Q}_{f,\Gamma}$. Les $S_{h}$ premières
singularités, pôles des néga-hémisphères, seront dite de type "h",
les $F$ dernières seront dite de type "f".
\end{Def}
On a maintenant l'analogue du théorème \ref{caracmet} dans le cas
fuchsien :
\begin{The}\label{caracmefuchs}
On a l'égalité $\mathcal{Q}_{f,\Gamma}=\mathcal{M}^{h}_{f,\Gamma}$.\\
Un espace métrique $h$ de cet ensemble avec ses $S_{h}+F$
singularités numérotées étant donnée, on sait retrouver
directement les géodésiques de $h$ formant la combinatoire duale
du polyèdre tronqué, image de $h$ par $\tilde{\phi}^{-1}$.
\end{The}
\paragraph{Démonstration.} Il faut vérifier que les étapes de la
démonstration du théorème \ref{caracmet} ne sont pas affectées par
le fait que l'on considère une combinatoire de genre $2$. On
remarque que tout le raisonnement de la démonstration du théorème
\ref{caracmet} utilise seulement deux choses :
\begin{enumerate}
\item le fait que l'espace métrique peut se construire comme un recollement d'hémisphères
\item la convexité du polyèdre.
\end{enumerate}
L'argument reste donc valable dans le cas fuchsien.$\diamond$\\
\subsection{Caractérisation de polyèdres fuchsiens hyperidéaux en fonction de leurs angles dièdres}
On peut maintenant caractériser les polyèdres fuchsiens en
fonction de leurs angles dièdres.
\begin{Def}\label{defKgammafuchs}
$K_{f,\Gamma}$ est le
sous-ensemble de $]0,\pi[^{E}$ défini de la manière suivante:\\
A chaque arête $e_{i}$ de $\Gamma^{*}$ on associe un poids
$\theta_{e_{i}}$ (représentant les angles dièdres) et vérifiant
les conditions suivantes :
\begin{itemize}
\item $\theta_{e_{i}} \in ]0,\pi[$
\item \emph{CF1}- $\sum_{i=1}^{n}\theta_{e_{i}}\geq 2\pi$ pour toute courbe
fermée contractile $\gamma$ plongée dans le graphe $\Gamma^{*}$ et
passant par $e_{1},...,e_{n}$; l'égalité ayant lieu très
exactement quand $\gamma$ est le bord d'une face du graphe dual
$\Gamma^{*}$ correspondant à un sommet idéal.
\item \emph{CF2}- $\sum_{i=1}^{n}\theta_{e_{i}} > \pi$ pour tout chemin $\gamma$ plongé
dans $\Gamma^{*}$ et joignant deux sommets $A$ et $B$ de
$\Gamma^{*}$ ayant une même face adjacente $F$ , mais tel que
$\gamma$ ne soit pas tout entier contenu dans le bord de cette
face. On demande que $\gamma$ soit homotope au chemin bordant $F$
et joignant $A$ à $B$.
\end{itemize}
\end{Def}
\begin{The}[Caractérisation des polyèdres hyperidéaux fuchsien.] \label{caracthypfuchs}L'application $\Theta$
qui à un polyèdre fuchsien de combinatoire $\Gamma$, associe ses
angles dièdres est un homéomorphisme de $\mathcal{P}_{f,\Gamma}$
vers $K_{f,\Gamma}$.
\end{The}
\paragraph{Démonstration.}Il faut et il suffit de vérifier que la
condition pour les métriques duales d'avoir des géodésiques
fermées contractiles de longueur strictement supérieur à $2\pi$,
sauf si elles bordent un hémisphère, est équivalente aux
conditions (CF1) et (CF2) sur les angles dièdres de la définition
\ref{defKgammafuchs}. Cette vérification est la même que dans le
cas non-fuchsien du théorème \ref{caracthyp}, en prenant soin de
traduire le fait que l'on a affaire aux géodésiques contractiles :
les chemins fermés dans le graphe de $\Gamma^{*}$ de la condition
(C1) doivent être contractiles, et les chemins de la condition
(C2) doivent être contractiles une fois rajouté un segment de
géodésique de l'hémisphère où aboutissent les
extrémités.$\diamond$
\section{Conclusion mathématique}
Lorsqu'on cherche à caractériser les polyèdres hyperboliques par
leurs angles dièdres, l'essentiel des problèmes vient du fait
suivant : La connaissance de la métrique (duale ou pas), ne permet
pas de reconstruire la combinatoire de l'unique polyèdre
correspondant à cette métrique, et ceci dés que l'un des sommets
est compact. Curieusement, on ne peut donc toujours pas par ce
moyen répondre à la question de rigidité dans le cas compact
général en fonction des angles dièdres, bien que l'on sache le
faire en fonction de la métrique duale.
\section{Remerciements}
Mes plus vifs remerciements vont évidemment à Jean-Marc Schlenker,
non seulement pour la très grande qualité du stage, qui m'a permis
de me familiariser modestement à la recherche, mais aussi pour la
qualité des relations humaines. C'est de lui que proviennent
beaucoup des idées qui ont donné lieu à cet article.

\end{document}